\documentclass[11pt]{amsart}
\usepackage{amsmath, latexsym, amssymb}
\input xypic
\begin{document}

\title[Monotone invariants and embeddings of  statistical manifolds]{ Monotone invariants and embeddings of  statistical manifolds}\footnote{appeared in ``Advances in  Deterministic  and Stochastic Analysis",  World Scientific Hackensack, NJ, 2007, 231-254.}
\author[H.V. L\^e]{  H\^ong V\^an L\^e \\
Max-Planck-Institute for Mathematics in Sciences\\
  Inselstra\ss e 22-26\\
  D-04103 Leipzig\\
emai: hvle@mis.mpg.de}
\maketitle
\newcommand{\R}{{\Bbb R}}
\newcommand{\C}{{\Bbb C}}
\newcommand{\F}{{\Bbb F}}
\newcommand{\Z}{{\Bbb Z}}
\newcommand{\N}{{\Bbb N}}
\newcommand{\Q}{{\Bbb Q}}

\newcommand{\Aa}{{\mathcal A}}
\newcommand{\Bb}{{\mathcal B}}
\newcommand{\Cc}{{\mathcal C}}    
\newcommand{\Dd}{{\mathcal D}}
\newcommand{\Ee}{{\mathcal E}}
\newcommand{\Ff}{{\mathcal F}}
\newcommand{\Gg}{{\mathcal G}}    
\newcommand{\Hh}{{\mathcal H}}
\newcommand{\Kk}{{\mathcal K}}
\newcommand{\Ii}{{\mathcal I}}
\newcommand{\Jj}{{\mathcal J}}
\newcommand{\Ll}{{\mathcal L}}    
\newcommand{\Mm}{{\mathcal M}}    
\newcommand{\Nn}{{\mathcal N}}
\newcommand{\Oo}{{\mathcal O}}
\newcommand{\Pp}{{\mathcal P}}
\newcommand{\Qq}{{\mathcal Q}}
\newcommand{\Rr}{{\mathcal R}}
\newcommand{\Ss}{{\mathcal S}}

\newcommand{\Uu}{{\mathcal U}}
\newcommand{\Vv}{{\mathcal V}}
\newcommand{\Ww}{{\mathcal W}}
\newcommand{\Xx}{{\mathcal X}}
\newcommand{\Yy}{{\mathcal Y}}
\newcommand{\Zz}{{\mathcal Z}}

\newcommand{\zt}{{\tilde z}}
\newcommand{\xt}{{\tilde x}}
\newcommand{\Ht}{\widetilde{H}}
\newcommand{\ut}{{\tilde u}}
\newcommand{\Mt}{{\widetilde M}}
\newcommand{\Llt}{{\widetilde{\mathcal L}}}
\newcommand{\yt}{{\tilde y}}
\newcommand{\vt}{{\tilde v}}
\newcommand{\Ppt}{{\widetilde{\mathcal P}}}

\newcommand{\Remark}{{\it Remark}}
\newcommand{\Proof}{{\it Proof}}
\newcommand{\ad}{{\rm ad}}
\newcommand{\Om}{{\Omega}}
\newcommand{\om}{{\omega}}
\newcommand{\eps}{{\varepsilon}}
\newcommand{\Di}{{\rm Diff}}
\newcommand{\Pro}[1]{\noindent {\bf Proposition #1}}

\newcommand{\Lem}[1]{\noindent {\bf Lemma #1 }}
\newcommand{\An}[1]{\noindent {\bf Anmerkung #1}}
\newcommand{\Kor}[1]{\noindent {\bf Korollar #1}}
\newcommand{\Satz}[1]{\noindent {\bf Satz #1}}

\newcommand{\gl}{{\frak gl}}
\newcommand{\g}{{mathfrak g}}
\renewcommand{\o}{{\frak o}}
\newcommand{\so}{{\frak so}}
\renewcommand{\u}{{\frak u}}
\newcommand{\su}{{\frak su}}
\newcommand{\ssl}{{\frak sl}}
\newcommand{\ssp}{{\frak sp}}

\newcommand{\Cinf}{C^{\infty}}
\newcommand{\CS}{{\mathcal{CS}}}
\newcommand{\YM}{{\mathcal{YM}}}
\newcommand{\Jreg}{{\mathcal J}_{\rm reg}}
\newcommand{\Hreg}{{\mathcal H}_{\rm reg}}
\newcommand{\SP}{{\rm SP}}
\newcommand{\im}{{\rm im}}

\newcommand{\inner}[2]{\langle #1, #2\rangle}    
\newcommand{\Inner}[2]{#1\cdot#2}
\def\NABLA#1{{\mathop{\nabla\kern-.5ex\lower1ex\hbox{$#1$}}}}
\def\Nabla#1{\nabla\kern-.5ex{}_#1}

\newcommand{\half}{\scriptstyle\frac{1}{2}}
\newcommand{\p}{{\partial}}
\newcommand{\notsub}{\not\subset}
\newcommand{\iI}{{I}}               
\newcommand{\bI}{{\partial I}}      
\newcommand{\LRA}{\Longrightarrow}
\newcommand{\LLA}{\Longleftarrow}
\newcommand{\lra}{\longrightarrow}
\newcommand{\LLR}{\Longleftrightarrow}
\newcommand{\lla}{\longleftarrow}
\newcommand{\INTO}{\hookrightarrow}
\newcommand{\la}{\langle}
\newcommand{\ra}{\rangle}

\newcommand{\Sy}{\text{ Diff }_{\om}}
\newcommand{\Ex}{\text{Diff }_{ex}}
\newcommand{\jdef}[1]{{\bf #1}}
\newcommand{\QED}{\hfill$\Box$\medskip}
\newcommand{\bm}{\boldmath}


\numberwithin{equation}{section}
\theoremstyle{plain}
\newtheorem{lemma}{Lemma}[section]
\newtheorem{proposition}[lemma]{Proposition}
\newtheorem{theorem}[lemma]{Theorem}
\newtheorem{corollary}[lemma]{Corollary}
\newtheorem{problem}[lemma]{Problem}

\theoremstyle{definition}
\newtheorem{definition}[lemma]{Definition}
\newtheorem{remark}[lemma]{Remark}
\newtheorem{example}[lemma]{Example}

\medskip

\begin{abstract}
In this note we  prove certain necessary and sufficient conditions
for the existence of  an embedding of statistical manifolds. In particular, we
prove that any  compact smooth ($C^1$ resp.) statistical manifold can be embedded
into the  space of probability measures on a finite set. As a
result, we get  an answer to the Lauritzen question on a realization of
smooth ($C^1$ resp.) statistical manifolds   as statistical models.

\medskip

\end{abstract}


{\it Keywords: Fisher metric, Chentsov-Amari connections, statistical\\ manifolds, statistical models.}

\medskip


\medskip

\section { Introduction}

A {\bf statistical model} is a family $M$ of probability measures on  a measurable space
$\Om$. There are two natural geometrical structures  on any statistical model equipped
with a differentiable manifold structure. They are the  Fisher  tensor and the Chentsov-Amari tensor.

The  Fisher tensor was given by  Fisher in 1925 as an
information characterization of a statistical model.  Rao [Rao(1945)] proposed to consider this
tensor as a Riemannian metric on the manifold of probability
distributions. This Fisher metric  has been systematically studied
in [Chentsov1972], [M-C 1990], [A-N2000] and others [Lauritzen1987], [Rao1987],  [Ay2002],
[Jost2005], ect. in the field of geometric aspects  of statistics
and information theory.

Chentsov [Chentsov1972] and Amari [Amari1997]  independently also
discovered a natural structure on statistical models, namely  a
1-parameter family of invariant  connections,  which includes the
Levi-Civita connection of the Fisher metric.  This family of
invariant connections is defined  by a 3-symmetric  tensor $T$
together with the Levi-Civita connection of the Fisher metric.

Motivated by the question how much we can describe a statistical model via
their Fisher metric and Chentsov-Amari tensor $T$, in 1987 Lauritzen proposed to
 call a Riemannian manifold $(M, g)$ with a
3-symmetric tensor $T$  {\bf a statistical manifold}. Since two
3-symmetric tensors $T$ and $k\cdot T$, $k\not=0$, define the same
family of Chentsov-Amari connections,  we shall say that two
statistical manifolds $(M, g, T)$ and $(M,g, kT)$ are conformal
equivalent.

A natural and important question  in the mathematical
statistics is to understand,
if a given family $M$ of probability distributions can be considered as a subfamily of
another given  one $N$.  In the language of statistical manifolds, this question can be formulated  as a problem of  isostatistical
embedding of a statistical manifold $(M, g, T)$
into  another one $(N, g', T')$. Here we say that an immersion $f: (M, g, T) \to (\bar M, \bar g, \bar T)$ is called
{\bf isostatistical}, if $ f ^* (\bar g) = g$ and $f^* (\bar T) = T$.

We shall see in  section 2 that the   problem  of the existence of an isostatistical embedding includes also  the Lauritzen  question in 1987,  if
any statistical manifold is a statistical model.
It also concerns the following important  problem posed  by
Amari in 1997,  if any finite dimensional statistical model   can be embedded into the space
$Cap^N$ of  probability distributions of the sample space $\Om^N$ of $N$ elementary  events
for some  finite $N$.

We shall  construct a class of  $C^0$ (and $C^1$) {\bf monotone invariants} of statistical manifolds, which
present obstructions to embedding of  a given  $C^k$ statistical manifold $M$ into another one
$N^n$. Here a $C^k$ statistical manifold $(M, g, T)$ is a smooth differentiable manifold
with $C^k$ sections $g \in  S^2 T^*M $ and $ T \in S^3 T^* M$.
These invariants measure
certain relations between the metric tensor $g$ and the 3-symmetric tensor $T$.
In particular,  using these invariants we show that no
statistical manifold which is conformal equivalent to the space $Cap ^N$  can be embedded
into  the product  of $m$ copies of the  normal
Gaussian  manifolds for any   $N >3$ and any  finite $m$. 
In the Main Theorem (section 5) we prove that any  compact smooth ($C^1$ resp.) statistical
manifold  $M ^m$ can be  isostatistically embedded to a the space $Cap ^N $ for  some
$N$ big enough.

As a consequence we  also get a new proof of Matumoto theorem on the existence of
the contrast function for a   compact statistical manifold (see 2.8).

\medskip

{\bf Acknowledgement}. I am thankful to J\"urgen Jost and Nihat Ay for
their introduction to the field of information geometry and  helpful discussions.

\medskip

\section{Statistical models and statistical manifolds.}

In this section we recall the definitions of
the Fisher metric and the Chentsov-Amari connections on statistical models.
We introduce the notion of a weak Fisher metric and a weak potential function.
At the end of the section we discuss the problem,
if a given statistical manifold is a statistical model. Most of the facts in this section
can be found in [A-N2000].

\medskip

Suppose that $M$ is a statistical model - a family of probability measures on
a space $\Om$. We assume throughout this note that
$M$ and $\Om$ are differentiable manifolds, and $\Om$ is equipped with a fixed Borel measure $d\om$. We  also
write 
$$p(x,\om) = p(x, \om) d\om,\leqno (2.1)$$
where $p(x,\om)$ in LHS of (2.1) is a Borel measure in $M$ and $p (x, \om)$ in the RHS of (2.1) is a non negative (density) function on $M\times \Om$ which satisfies
$$\int_{\Om} p (x,\om) \, d\om = 1\: \forall x\in M. \leqno (2.1.a)$$
The {\bf Fisher metric} $g^F (x)$ 
is defined on $M$
 as follows. For any $V, W \in T_x M$ we put
$$ g^F ( V, W) _x =  \int _{\Om}(\p_V \ln p(x, \om))(\p_W \ln p(x, \om)) p(x, \om).\leqno(2.2)$$

 The function under integral
in (2.2)
is well defined, if 
$$ p(x,\om) > 0  ,\leqno (2.1.b)$$

Denote by $Cap (\Om)$ the space of all probability measures on $\Om$.
Clearly we can consider the  density function $p (x,\om)$ as  a mapping $M\to Cap (\Om)$. Thus we shall  call a function $p(x, \om)$
{\bf a probability potential of the metric $g_F$}, if $p(x,\om)$ satisfies (2.1.a), (2.1.b), (2.2).  (It is known  that
for a given   Riemannian metric $g_F$  on a smooth manifold $M$  there  exist
many probability potentials $f(x, \om)$ for $g_F$, even if we fix  the space
$(\Om,d\om)$.)

Some time it is useful to consider functions $p(x,\om)$ which satisfy (2.2) and (2.1.b) but not
necessary (2.1.a).  In this case, the  Riemannian metric $g ^F$ will be called {\bf weak
Fisher metric}, and  the function $p(x,\om)$ will be called {\bf a weak probability potential}
of $g^F$.

\medskip

{\bf 2.3. Example of a weak Fisher metric: the  standard Euclidean metric $g^0$ } on the positive quadrant
 $\R ^N _+ ( x_i> 0)$. It is straightforward to check that $g_0$ admits a weak probability potential
$\{ p_i (x) = {1\over 4} x_i ^2, i = \overline{1,N}.\}$ 
Here $\Om = \Om^N$ - the sample space  of $N$ elementary events.
\medskip

{\bf 2.4. The Fisher metric on the space {$ \bm (Cap  ^N)_+ $}}  of all  positive probability
distributions on $\Om^N$ (see also [A-N2000], [Jost2005], [Chentsov1972]).
By definition we have
$$ Cap ^N_+ : = \{ ( p_1, \cdots, p_N)| p_i  > 0 \text{ for }
i=\overline{1,N} \: \& \sum p_i = 1\}.$$
We  define the  embedding map
$$ f : Cap ^N_+ \to S^{N-1}(2),$$
$$( p_1, \cdots, p_N)\mapsto  (q_1=2\, \sqrt{p_1}, \cdots,  q_N= 2\, \sqrt{p_N}).$$
It is easy to see that the Fisher metric  in the new coordinates $(q_i)$ is the standard metric of  constant positive curvature on the sphere $S^{N-1}(2)$.

\medskip

{\bf 2.5. Divergence  potential} (see [A-N2000],  [Rao(1987)].)   A  function $\rho$ on $M\times M$ with the following property
$$\rho (x,y)g^Fe 0 \text { with  equality  iff }  x = y\leqno (2.5.1)$$
is called a divergence function.
A divergence function $\rho$  is called {\bf a divergence potential} for
a metric $g$ on $M$, if 
$$ g (X, Y) _ x = Hess (\rho) ( i_1 (X), i_1 (Y)).\leqno (2.5.2)$$
where
$$T_{(x, x)} (M, M) = ( T_xM, 0) \oplus (0, T_xM) = ( i_1 (T_xM)) \oplus (i_2(T_xM)).$$
 An example of a
divergence potential  for a Fisher metric is  the Jensen function  $J^{\lambda, \mu}_H (x,y)$
of the entropy  function  $H(x)$ on $M$,
or a Kullback relative entropy function $K(x,y)$ on $M \times M$.  


{\bf 2.6. Chentsov-Amari connections.}  Let $p(x,\om)$ be a  probability potential  for a Riemannian
metric $g$. We define a symmetric 3-tensor $T$ on $M$ as follows
 $$ T(X, Y, Z) = \int (\p _X \ln p (x, \om)) (\p _Y\ln  p ( x, \om) )(\p _Z \ln p (x, \om) )\,
 p (x, \om).\leqno(2.6.1)$$
We denote by $\nabla ^F$ the Levi-Civita connection of the (weak) Fisher  metric  $g^F$.
We define
$$<\nabla ^t_X Y, Z > : = <\nabla ^F _ X Y, Z > + t\cdot T (X, Y, Z).\leqno(2.6.2)$$
The connections $\nabla^t$  are called the Chentsov-Amari connections.

\medskip

{\bf 2.6.3.  Remark.} ([A-N2000], [Matsumoto1993]) Any divergence function $\rho(x,y)$ on $M\times M$ defines a  tensor  $T$ on $M$ via the following formula
$$ T(X, Y, Z) _{x} =- \p _{ i_2 (Z)} Hess (\rho ) ( i_1(X), i_1 (Y))_{(x,x)}+ \p _{i_1 (Z)}
Hess(\rho)  (i_2(X), i_2(Y))_{(x,x)}.$$
If $g$ and $T$ are defined by
the same divergence
function $\rho (x,y)$, we  shall call $\rho (x,y)$
{\bf a divergence potential for  the  statistical manifold $(M,g,T)$}.
It is a known fact that  the Kullback relative entropy function is a divergence potential
for the associated statistical model. 

\medskip

{\bf  2.7. Statistical submanifolds.} 

\medskip

A submanifold $N$ in a statistical manifold $(M, g, T)$ with
the induced  Riemannian metric $g_{| N}$ and induced
tensor $T_{| N}$ is called {\bf statistical submanifold} of
$(M, g, T)$. Clearly, if $f(x, \om)$ is a  (weak) probability potential for $(M, g, T)$,
then its restriction to any submanifold $N\subset M$ is a (weak) probability
potential of the induced statistical structure.

\medskip

{\bf 2.8.  Statistical models and statistical manifolds.}
 Since any probability function  $p(x,\om)$ defines a map $M\to Cap (\Om)$, we shall
say that a statistical manifold $(M,g,T)$ is a statistical model, if there
probability potential
$p(x,\om)$ for $g$ and $T$. By the remark in 2.7, we get that a statistical submanifold of a statistical model is also a statistical model. Furthermore,  if a statistical manifold $(M,g,T)$ is a  statistical
model, then it must admit a divergence potential. Hence   we obtain the following

\medskip

{\bf 2.8.1  Theorem. } (cf. [Matumoto1993] ) {\it  For any  compact statistical manifold $(M,g, T)$ there exists
a  divergence potential $\rho$ for $g$ and  for $T$.  }

Note that  Matumoto's theorem  does not requires the compactness  of  $(M, g, T)$.

\medskip

\section{ Embeddings of linear statistical spaces.}

An Euclidean space $ (\R^n, g ^0)$ equipped with a 3 -symmetric  tensor
$T$ will be called {\bf a linear statistical spaces.} We observe that
the equivalence class of   linear statistical spaces coincides with the orbit
space of 3-symmetric tensors $T$ under the action of
the orthogonal group $O(n)$.  In this section
we discuss certain invariants of these orbits and we show several necessary and
sufficient conditions for the existence of  embedding of one linear statistical space into
another linear statistical space by studying these invariants. A class of our necessary conditions consists
of {\bf monotone invariants $\lambda$}, i.e. we assign to any linear statistical  space
$(\R^n, g ^0, T)$ a number  $\lambda (\R^n, g^0, T)$ such that, if $(\R^n, g^0, T)$
is a statistical submanifold of $(\R^m, g^0, T')$, then we have
$$\lambda (\R^n, g^0, T) \le \lambda (\R^m, g^0, T').$$
 Since a tangent space of a statistical manifold is a linear statistical manifold, these invariants play important role in the problem of isostatistical immersion.
\medskip

{\bf 3.1. Trace type of a symmetric 3-tensor.} 
Let us denote by $\Rr ^n$  the subspace in $S^3(\R^n)$ consisting of the following
3-symmetric tensors
$$ T^v (x, y, z) = <v, x><y,z> + <v, y><x,z> + <v,z><x,y>,$$
where $v \in \R ^n$.
Using the standard representation theory (see e.g. [O-N1988]) we have the decomposition
$$S^3(\R ^n)  = \Rr (3\pi_1)\oplus \Rr ^n.\leqno (3.2)$$
\medskip

The  component $\Rr^n$ is defined by taking the trace of $T$
$$Tr: S^3 (\R^n)^* \to (\R^n)^*, \: Tr (T) (v) : = Tr (v \rfloor T).$$
Clearly $Tr$ is an $SO(n)$-equivariant map with nonzero image. Using the identity   $Tr ( T^v) = (n+2) v^*$,  we get

{\bf  3.3.   Lemma.} {\it We have}
$$\pi_2 (T) = {1\over n+2} T ^{ Tr (T)}.\leqno(3.4)$$
\medskip

In view of Lemma 3.3 we shall call any tensor $T\in\Rr^n$ of {\bf trace type.}

\medskip

We note that
$$\dim S^3 ( \R^n) = C^3_n + 2 C^2_n  + n = {n(n+1)(n+2)\over 6}.$$
Thus  the dimension of the quotient $S^3 (\R^n) / SO (n)$ is at least
$C^3_n + C^2_n + n$. A direct computation shows that the dimension
of the orbit $SO(n) ( [ \sum _{i=1} ^n a_i v_i ^3])$ is $ C^2_n =\dim SO(n)$, if
$\Pi a_i\not= 0$.
Here $\{v_i\}$ is an orthonormal basis in $\R^n$. Hence the
dimension of $S^3(\R^n)/ O(n) = C^3 _n + C^2_n + n$. This dimension
is exactly the number of all complete invariants of pairs consisting  of a positive definite
bilinear form $g$ and a 3-symmetric  tensor $T$.

\medskip

Since the dimension of $G_k (\R^n) = k(n-k)$, it follows that  generically
it is impossible to embed a linear  statistical space $ (R^k, g^0, T)$ into
a  given statistical  linear space $(R^n, g ^0, T)$, unless $ k(n-k) g^Fe
C^3_k + C^2_k + k$.  Clearly the dimension condition
is not sufficient as the following proposition shows.

\medskip

{\bf 3.5.  Proposition.} {\it A linear statistical space $(\R^k, g^0, T)$ can be
embedded into a linear statistical space $( \R^N, g^0, T ^v)$, if and only
if $Ng^Fe k$ and  $T$ is also a trace type: $T = T ^w$ with $|w| \le  | v|$.}

\medskip

{\it Proof.} The necessary condition follows from the fact that the restriction of $T^v$ to $\R^k$  equals
$ T^{\bar v}$, where $\bar v$ is the  orthogonal projection of  $v$ to $\R^k$.
Conversely, if $|w| \le |v|$ we can find an orthogonal transformation, such
that $w$ equals the orthogonal projection of $v$ on $\R^k$.\QED

\medskip

 {\bf 3.6. Commasses as monotone invariants.} Since the metric $g$  extends canonically
on the space $S^3 (\R^n)$, we can define the absolute norm
$$ || T ||: = \sqrt {<T, T >}.$$
 Now we define   {\bf comasses}  of a 3-symmetric tensor $T$ as follows
$$\Mm^3 ( T) : = \max_{ |x|= 1,|y| = 1,  |z| = 1} T (x, y, z),$$
$$ \Mm ^{2} (T)  : = \max _{|x| = 1,    | y| = 1} T (x, y, y),$$
$$\Mm ^1 (T) : = \max _{| x| = 1} T (x,x,x).$$
 Clearly we have
$$ 0\le \Mm^1(T)\le \Mm ^2 (T) \le \Mm ^3(T) \le || T|| .$$
\medskip

{\bf 3.7. Proposition.} {\it The   comasses   $\Mm^i$, $ i \in [1, 3]$,  are  nonnegative  linear monotone invariants, which vanish  if and only if $T = 0$.}

\medskip

{\it Proof. }  Clearly $\Mm^i(T) \ge 0$ for $i = 1,2,3$.
Now we are going  to show that
$\Mm^1$ vanishes at  $T$ only if $ T = 0$. Observe  that  $\Mm^1 = 0$ if and only if $T(x, x, x) = 0$ for  all $x \in \R^n$.  Writing  $T$ in  coordinate  expression  $T (x, y, z) = \sum a_{ijk} x^i y ^j z ^k$, we   note that
$T(x, x, x) = 0$ if and only if $T = 0$, since  $T$ is symmetric.

 Next   we shall show  that  $\Mm^i(T)$  is  a linear   monotone invariant   for $i=1,2,3$.
Assume that $e $ is a linear  embedding $(\R^n, g, T)$   into $(\R^m, \bar g, \bar T)$.
Then $T$ is a restriction of the 3-symmetric tensor $\bar T$.  Hence  we have
$$  \Mm^i (T) \le \Mm ^i ( \bar T), \, \text{ for } i = 1,2,3.$$
This  implies  that    $\Mm^i$  are linear  monotone invariants.
\QED
\medskip

Now for a space $(\R^n, g^0, T)$  and for $ 1 \le k\le n $ we put
$$\lambda _k (T) := \min_{\R^k\subset \R ^n} \Mm ^1( T_{| \R ^k}).$$

We can easily check  that  if $\bar T$ is a restriction of $T$ to a subspace
$\R ^m\subset \R^n$, then
$$\lambda _k (\bar T) \ge \lambda _k(T) \ge 0 \text { for all }  k\le m.$$

Thus $-\lambda_k (T)$ is a monotone invariant of   linear statistical manifolds. These
invariants are related by the following inequalities
$$\Mm^1(T) = \lambda _n (T) \ge \lambda_{n-1} (T) \cdots \ge \lambda _2 (T) \ge
\lambda_1 (T) = 0.$$

The last equality follows from the fact, that the function $T(x,x,x)$ is anti-symmetric
on $S^{n-1}(|x| = 1)\subset\R^n$ and $S^{n-1}$ is connected. We observe that if $T$ is
of trace
type, then $\lambda_{n-1} (T) =\cdots  = \lambda_1 (T)  = 0$.

\medskip

We  are going to give a lower bound of the monotone invariant $\lambda_{n-1}$ of
a
linear statistical space of  certain type.
The equality $\lambda_{n-1}(\R^n, g^0, T)  g^Fe A$ means that no hyperplane with the norm $\Mm ^1$
strictly less  than $A$ can be embedded in $(\R^{n-1}, g^0, T)$.

\medskip

{\bf 3.8. Lemma.} {\it a) Let  $T = \sum_{i = 1 } ^n (N -\eps _i ) (x ^i) ^3$ be a 3-symmetric
tensor on $\R ^n$  with $n \ge 4$, $N \ge 4$ and $|\eps _i| \le 1 /4$.   Then  we have
$$\lambda_{n-1}( T) \ge  {N\over \sqrt 10}- 1/4.$$
b ) Let  $ T = N\sum _{ i=1} ^n  (x ^ i) ^ 3$, and  $H$  be a hyperplane in
$\R ^n$ which is orthogonal  to $(kn, 1,1, \cdots, 1)$, and  let $n \ge 5, \, k \ge 3$.
 Then  we have
$$\lambda _{ n-2} ( T _{ | H} )
\ge { N \over 5} -1.$$
c) Let $ x = ((1-\eps),  {1\over kn}, \cdots , {1 \over kn})\in S^n (1) \subset \R^{n+1}$,
where $n\ge 4, k  \ge n$.  We denote by $H$ the tangential plane $T_x S^n$, and
 by $T^0$ the  following   3-symmetric tensor on $\R^{n+1}_+$:
$$T^0_{ijk} (x_1, \cdots x_N) = \delta _{ijk} {2\over x_i}.\leqno (3.8.1)$$
Then we have}
$$\lambda _{n-1} ( T^0_{| H}) \ge {kn\over 5} -1.$$

\medskip

{\bf 3.8.2. Remark.}  The tensor $T^0$ in (3.8.1)  defines  on $(\R^n, g^0)$
a statistical structure with a weak probability potential $\{ {1\over 4} x_i^2,  i =1,n\}$, see  also 2.3.

\medskip

{\it Proof of Lemma 3.8.} The reader shall see that  a proof of Lemma 3.8 can be done in the same scheme of  the proof of Sublemma 5.10.
Therefore we do not repeat this argument here.

\medskip

{\bf 3.8.3. Remark.} {\it Lemma 3.8.a holds also for $n =3$ but not for
$n=2$, Lemma 3.8.b
holds also for  $n =4$, but not for $n=3$, and Lemma 3.8.c holds also for $n =3$
but not for $n=2$.}

\medskip

There are  also several  obvious  monotone invariants of $T$.
$$A^1( T) : = \max_{|x|= |y|  = |z| =1, <x,y> = <y,z> = <z, x>  = 0} T(x,y,z)$$
is well-defined for $n\ge 3$.
$$ A^2 (T) := \max _{|x| = |y| = 1, <x,y > = 0} T (x, y, y),$$
is well-defined for $n\ge 2$.
We can check that
$$ \ker  A^1 = \Rr ^n .$$
On the other hand we have
$$\ker A^2 \subset \Rr ( 3\pi_1).$$

Thus $A ^1$ and $A^2$ are different invariants.

\medskip

{\bf  3.9. Lemma.} {\it Let $\pi_1$  be the first component  of  $T$ in decomposition (3.2).
Then
$|| T|| _1 : = ||\pi _1 (T)||$ is  a  monotone invariant of $T$.}
\medskip

{\it Proof. }  Let $\R^k$ be a subspace of $\R^n$. We denote
by $\pi_k ^n T$ the restriction of $T$ to $\R^k$. Clearly
$$\pi_k ^n (T) = \pi_k ^n ( \pi_1 T) + \pi_k ^n (\pi_2 T).$$
We have noticed in Proposition 3.5  that the restriction of  the trace form $\pi_2 T$ to any
subspace is also  a trace form. Thus $\pi_k ^n (\pi_2 )$ is an element in $ \Rr ^k \subset
S^3 (\R ^k).$
Hence we have
$$\pi_1 ( \pi_k ^n T) = \pi _1 (\pi_k^n ( \pi_1 T)).\leqno (3.9.1)$$
Since all the projections $\pi_1$, $\pi_k ^n$ decrease the norm $||. ||$, we get
$$|| \pi_k ^n T||_1 = || \pi_1 (\pi_k ^n T)||  = ||\pi_1 (\pi_k ^n(\pi_1 T))||
 \le ||\pi_ 1 (T) || = ||T|| _1.$$
\QED

{\bf 3.10. Proposition.} {\it   A  statistical line
$ (\R, g^0, T)$  can be embedded into $( \R^N, g^0, T ')$, if and only if $\Mm ^1(T)\le \Mm^1 (T')$}.

\medskip
{\it Proof.}  It suffices to show that we can embed $(\R, g^0, T)$ into $(\R^N, g^0, T')$, if
we have $\Mm^1(T) \le \Mm^1 ( T')$.   We note that $T'(v,v,v)$ defines an anti-symmetric
function on the sphere $S^{N-1}(|v| = 1)\subset  \R^N$. Thus there is a point $v\in S^{N-1}$ such that
$T' (v,v,v) = \Mm^1 (T)$. Clearly the line $v\otimes \R $ defines the required embedding.\QED

\medskip

 Let us consider  the embedding problem for
2-dimensional linear statistical spaces. It is easy to see that
$$S^3 (\R^2) = \R^2 \oplus \R^2.$$
Thus the quotient $S^3 (\R^2)/ SO(2)$ equals $( \R^2 \oplus \R^2)/ S^1$.  Geometrically there
are several ways to see this.
In the first way we denote components of $T\in S^3 (\R^2)$ via
$T_{111}, T_{112}, T_{122}, T_{222}.$
\medskip

{\bf 3.11. Lemma.} {\it There exists an oriented orthonormal basic in $\R^2$ such that
$T_{111}  = \Mm^1(T) > 0, T_{112} = 0$ for all non-vanishing $T$. These numbers
$(T_{111}, T_{122}, T_{222})$ are called
canonical coordinates of $T$. Two tensors $T$ and $T'$ are equivalent, if and only if they
have the same canonical coordinates.}
\medskip

{\it Proof.}
We choose an oriented orthonormal basis $(v_1, v_2)$ by taking as $v_1$  a
point on $S ^1(|x| =1)$, where the function $T(x,x,x)$ reaches the maximum. The first variation
formula shows that in this case $T_{112}=0$. This shows the existence of  the canonical coordinates.
Clearly, if two tensors have the same canonical coordinates, then they are equivalent.
Next, if two tensors $T$ and $T'$ are equivalent,  then their norms $\Mm^1$ are the same. We need to
take care the case,  when there are several points $x$ at which  $T(x,x,x)$ reaches the maximum. In any case,
they have the same first coordinates.  Next we note that
$$<Tr (T), Tr (T) > = (T_{111}+T_{122})^2 + T_{222} ^2,$$
$$ || T || ^2 = T_{111} ^2 + T_{122} ^2 + T_{222} ^2.$$
Thus if two tensors are equivalent and have the same first coordinates, they must
have the same third coordinate $T_{122}$, and this third coordinate is uniquely defined up to sign. The
condition on the orientation  tells us that the sign must be $+$. This proves
the second statement.\QED

\medskip


\medskip

{\bf 3.12. Proposition.}  {\it We can always embed the 2-dimensional
statistical space $( \R^2, g^0, 0)$ into any  linear statistical space $ (\R^n, g^0, T)$,
if $ n \ge 7$.}

\medskip

{\it Proof.}   It suffices to prove for
$n = 7$. We denote by $\Oo (T)$ the set of of all unit vectors $v\in S^6$
such that $T (v,v,v) = 0$. Clearly $\Oo(T)$ is a set of dimension 5 in
$S^{6}$. Since $T$ is anti-symmetric, there exists a connected
component $\Oo ^0 (T)$ of $\Oo (T)$ which is invariant under the anti-symmetry
involution. Now  we consider the following   function $f$ on $\Oo^0 (T)$.
For each $v \in \Oo ^0 (T)$ we denote by $A ^v$ the bilinear symmetric 2-form
on the space  $T_x \Oo ^0 (T)$ considered as a subspace in $\R^n$:
$$ A^v ( y,z) = T (v, y,z).$$
Then we define $f(v)$ equal to $ \det (A ^v)$.  Since  $\Oo (T)$  has dimension
5, the function $f(v)$ is anti-symmetric on $\Oo  ^0 (T)$. Hence the set $\Oo ^0 _0 ( T)$
of all $v\in \Oo ^0 (T)$ with $ f(v)=0$ has dimension 4 and  it
contains a connected component which is also  invariant under the anti-symmetric
involution. For the simplicity we denote this connected component
also by $\Oo ^0 _0 (T)$.   Now we consider the following  two  possible cases.

\underline{Case 1}. We  assume that there is a point $v\in \Oo ^0 _0(T)$ such that
the nullity of $A^v$ is at least 2. Then there are two linear independent vectors
$ y, z \in T_v $ such that  the  restriction of $A^v$ on the plane $\R^2(y, z)$
vanishes.  Since the set $\Oo ^0(T)$ is connected and anti-symmetric and of
co dimension 1 in $S^{n-1}$,  the plane
$\R (y,z)$ has a non-empty intersection with  $\Oo ^0(T)$ at a point $w$. Then
the restriction of $T$ on the plane $\R^2 ( v, w)$  is  vanished, because
$$ T (v,v,v) = T ( w, w, w) = 0$$
$$ T( v,w, w) = 0 \, ( \text {since  } A^v ( w,w) = 0),$$
$$ T ( v, v, w) = 0 \, ( \text {since } w \in T_v \Oo ^0 (T)).$$

\underline{Case 2}. We assume that  the nullity of $A ^v$ on $\Oo ^0_0(T)$ is constantly 1.
Using the anti-symmetric property of $A^v$  we conclude that the restriction
of $A ^v$ to the  plane $\R^4 (v)$ which is orthogonal to the  kernel
of $A^v$ has index  constantly 2. Thus there exists a vector $z$ which
is orthogonal to the kernel $y$ of $A ^v$ such that $A^v ( z,z) = 0$.
Clearly the restriction of $A^v$ to the plane $\R^2 (y,z)$ vanishes.
Now we can repeat
the argument in the case 1 to get a vector $w$ such that the restriction
of $T$ to $\R^2(v,w)$ vanishes.\QED

\medskip

{\bf   3.13.  Theorem.} {\it a) Any statistical
space $(\R^n, g^0, T)$ can be embedded in the statistical space
$( \R ^{ n(n+1)}, g^0,  T ' = 2  ||T||\sum _{i=1}
^{N(n)} x_i ^ 3)$, where $x_i$ are  the canonical Euclidean coordinates on $\R^{n(n+1)}$.\\
b) The trivial space $ ( \R ^ n, g ^0, 0)$ can be embedded into $ (\R ^{2n}, g ^0,
\sum _{ i=1} ^ {2n} ( dx ^ i) ^ 3)$ for all $n$.}

\medskip

{\it Proof}. a)
We prove by induction. The statement for $n = 1$ follows from Proposition 3.8.
Suppose that the statement is valid for all $n \le k$.

\medskip

{\bf  3.14. Lemma.} {\it  Suppose that $T \in S^3( \R^{k+1})$. Then there are orthonormal coordinates
$x_1, \cdots, x_{k}$ such that}
$$  T =x _1  \sum _{i =1} ^{k+1} a_i x_i ^2  + \sum_{1 < i,j,k} a_{ijk} x_i x_j x_k.
\leqno (3.14.1)$$

\medskip

{\it Proof of Lemma 3.14}.  We choose $v_1$  as the unit vector in $ S^k \subset\R^{k+1}$, on which the function
$T(v,v,v)$ reaches the maximum on the unit sphere $S^k$. The first variation formula shows that
$T(v_1, v_1, w) = 0$ for all $w$ which is orthogonal to $v_1$.  We denote by $\R^k$ the orthogonal complement
to $\R \cdot v_1$. Now we consider a bilinear symmetric
form  $A$ on $\R^k$ defined as follows
$$ A(x, y) = S(v_1, x,y).$$
There is an orthonormal basis on $\R^k$, where we can write $A(x, y) =
\sum _{i= 2} ^{k+1} a_i x_i ^2$. Clearly in this orthonormal basis we can write
$T$ in the form  in (3.14.1).\QED

\medskip

 {\it Continuation of the proof of  Theorem 3.13.a } We shall show  explicitly that  that any  statistical
space
$(\R^2, g^0, T = a_2 x_1 (x_2) ^2)$ can  be embedded in $(\R^4, g^0,
\sum _{i=1} ^4 (y_i) ^3)$, if $0 \le |a_2| \le 1 / 2$.
 We put
$$ L (v_1): =\pm ( {1\over  2}, {1\over  2}, -{1\over  2} ,-{1\over 2} ) \leqno (3.15.1)$$
$$ L(v_2 ):=( \sqrt {{1+ 2a_2 \over 2}}, - \sqrt{{1+2a_2\over 2}} ,\sqrt{{1-2a_2\over 2}}
 ,- \sqrt {{1-2a_2\over 2}}  ).\leqno(3.15.2)$$
Here we take the sign $+$ in (3.16.1), if $a_2 > 0$, and we take the sign $-$, if $ a_2 < 0$.
Clearly, $L$ defines the required embedding $\R ^2 \to \R ^ 4$.

This together with Proposition 3.8 and the induction assumption completes the proof of
Theorem 3.13. a.

\medskip

{\it Proof of Theorem 3.13. b}. We  decompose the embedding $ f :   ( \R ^n , g ^0, 0)$
to $ (\R ^{2n}, g ^0, \sum _{ i=1} ^{2n}(x^i)^3)$ as  follows
$$ f( x_1, \cdots, x_n) = ( f ^ 1( x_1), \cdots , f ^n ( x_n))$$
where $ f ^ i $ embeds the line $ (\R , ( dx ^ i) ^2, 0)$ into $( \R ^ 2 , ( dx  ^ { 2i -1})  ^ 2 + (dx ^{2i} ) ^ 2,
(dx ^ { 2i -1} ) ^ 3 + ( dx ^{ 2i} ) ^ 3)$. Clearly, $f$  is the required embedding.\QED

\medskip


\medskip

\section {Monotone invariants and obstructions to embeddings of statistical manifolds}


 Let $ K ( M, e)$ denote the  category of statistical manifolds   $M$ with morphisms being
embeddings.  Functors of  this category are called
{\bf monotone invariants} of statistical manifolds.  Clearly any monotone invariant is an  invariant of
statistical manifolds.

\medskip

{\bf 4.1. Examples}. There are many  monotone invariants which arise  from our
analysis in section 3.

a) {\bf Trace type of a statistical manifold.} A statistical manifold $(M, g, T)$ will be called of {\bf trace type}, if for all $x\in M$ the form $T(x)$ is of trace type (see 3.1.) It follows from
Proposition 3.5 that   any statistical submanifold  of a statistical
manifold of  trace type is also of trace type. Thus the trace type is
a monotone invariant. In particular we cannot embed the
statistical space $Cap ^N$
and the normal Gaussian space into any  statistical space of trace type.
On the other hand, unlike
the linear case, we cannot embed a statistical manifold  of trace type into another
one of trace type,
even if the norm condition is satisfied. For example, if the trace form is closed (or  exact),
then the trace form of its submanifolds is also closed  (resp. exact). Hence within
a class of statistical manifolds of trace type we get a new monotone invariants which
can be expressed via the closedness and the cohomology class of the corresponding trace form.

\medskip

b) {\bf Decomposability } of a statistical manifold.
We note that the class of
3-symmetric tensors of trace
form is a subclass of all {\bf decomposable  tensors} $T^3$ which are a symmetric  product of
1-forms
and  symmetric 2-forms. Any statistical submanifold of a statistical manifold with
a decomposable tensor $T$ has also the (induced) decomposable  tensor. Thus the decomposability
is also a monotone invariant. The Gaussian normal
2-dimensional manifold is an example of decomposable type but not of trace type.

\medskip

c) {\bf Rank and comass.} We define for any statistical manifold $(M, g, T)$ the following  number
$$rank ( T) = \sup  rank (T(x))$$
$$||T|| _0  = \sup _{x\in M} || T(x)|| .$$
$$ \Mm^1 (T) _0 = \sup_{x\in M} \Mm^1 (T(x)).$$
$$|| T||_{1,0} = \sup_{x\in M} || T(x) ||_1.$$

Clearly  these four numbers are monotone invariants of statistical manifolds.

\medskip

We recall that the normal Gaussian statistical manifold is the two dimensional
statistical model  which is upper half of the plane $\R^2 (\mu, \sigma)$ with the potential
$$ p(\mu, \sigma)(x) = { 1\over \sqrt{ 2\pi}\, \sigma} \exp ( { - ( x-\mu) ^2\over 2 \sigma ^ 2} ), $$
here $x \in \R$.

{\bf 4.2. Proposition.} {\it Any statistical manifold which  is conformal equivalent to the space
$Cap ^N$ cannot be  embedded into the direct product of $m$ copies of the normal Gaussian
statistical manifold 2.3.3.a for any $N\ge 3$ and  finite $m$.}

\medskip

{\it Proof.} It is easy to check  that $\Mm^1 ( Cap ^N) = \infty$.
Thus any statistical manifold which is conformal equivalent to $Cap^N$ has also
the infinite invariant $\Mm^1$.
On the other hand,  we compute easily  that
the norm $\Mm^1$ of the Gaussian normal manifold, as well as the norm $\Mm ^ 1$ of  a direct product
of its finite copies, is finite. Namely the norm $\Mm^1(\mu, \sigma)$
is  $\sqrt 2$ for all $(\mu, \sigma)$.
\QED

\medskip

{\bf 4.3. Diameters of  statistical manifolds.} For  a positive number $\rho > 0$  and a
statistical manifold $(M, g, T)$ we set
$$ d _\rho ( M, g , T): = \sup \{l\in R^+\cup \infty\, | \, \exists \text {  an  immersion
of } ([0, l], dx ^2, \rho(dx ) ^3) \text  { to } (M,g, T).\} $$
We shall call $d_\rho (M,g, T)$ {\bf the diameter with weight $\rho$} of $(M, g, T)$.
Clearly $d_\rho$ are monotone invariants for all $\rho$.

\medskip

 To estimate  the  diameter with weight $\rho$ of a given statistical manifold
$(M, g, T)$ we can proceed as follows. For each point $x \in M$ we denote by $D_{\rho}(x)$
the set of all unit tangential vector $v\in T_x M$ such that $T (v,v,v) =\rho$. We denote
by $D_{\rho} ^i (x)$ the connected components of $D _{\rho} ^i (x)$.
We say that  a unite vector
$v$ in $T_x M$
is {\bf $\rho$-characteristic} with weight $c(x)$, if  there exists $i $ such that
we have
$$c( x) = \min_{w \in D ^i_{\rho}(x)} < v, w> > 0.$$


We shall say that a point $x \in M$ is { \bf $\rho$-regular}, if there is an open neighborhood
$U_\eps (x) \subset M$ such that $ D_{\rho}( U_\eps) = U_\eps \times D_\rho (x)$.
It is easy to see that the set of all $\rho$-regular points is open and dense in $M$
for any given $\rho$.

\medskip

{\bf 4.4. Proposition.} {\it The  diameter $d_{\rho}$ of
$ (M^m, g, T)$ is infinite, if $m\ge 3$ and there exists a number $\eps >0$
such that    one of the following 2 conditions holds:\\
a)There exists  a $(\rho+\eps)$-regular point $x\in M$ such that
the convex hull $Cov ( D_{\rho+\eps} ^i (x)) $ of one of connected
components $D_{\rho+\eps}  ^ i (x)$ contains the origin point $ 0 \in T_x M^m$ as it interior point.\\
b) $(M^m, g, T)$ has   a complete Riemannian submanifold $(N,\bar g)$ such that
there exists a smooth section
$ x\mapsto (D_{\rho +\eps} (x)\cap TN)$
over $N$.}

\medskip

{\it Proof}. The  statement  under the first condition a) is based on the fundamental
Lemma of the
convex integration technique of Gromov.  Namely Gromov  proved that
[2.4.1.A, Gromov(1986)], if the convex hull of some path connected subset $A_0 \subset
\R^q$ contains  a small neighborhood of the origin,  then there exists a map
$f : S^1 \to \R^ q$ whose derivative sends $S^1$ into $A_0$.

\medskip

{\bf 4.5. Lemma.} {\it Under the condition
in Proposition 4.4.1 there exists a small neighborhood $U_\delta (x)$ in
$M$  and an embedded  oriented curve $S^1\subset U_\delta (x)$  such that
for all point $s(t) \in S^1$ we have $\Mm^1 (T_{s(t)} S^1) \ge \rho + (\eps/2)$.}

\medskip

{\it  Proof of Lemma 4.5.}  We denote by $Exp$ the exponential map
$T_x M^m \to M^m$ and by $DExp$ the differential of this exponential map
restricted to $S^{m-1}\times T_x M^m \subset T(T_x M^m)$.  Here
 $S^{m-1}$ is  the unit sphere in $T_x M^m$. The
space $T_xM^m$ is a linear statistical space, so we denote
by $\Mm^1_x$ the induced norm-function on $S^{m-1}\times T_xM^m$ as follows:
$$\Mm^1 _x (l) = T _x(l,l,l).$$
Since $DExp$ is a continuous function, whose  restriction
to $S^{m-1} \times 0$ is the identity, there exists   a ball $B(0,\delta)$
with center in $0\in T_x M$  such that
$$\Mm^1( DExp (l)) - \Mm^1_x (l) ) < \eps / 4 \leqno(4.5.1)$$
for all $l \in S^{m-1} \times B(\delta) \subset T(T_x M^m)$.
We can assume that $\delta$ is so small such that $DExp$ is a homeomorphism
on $S^{m-1}\times B(0,\delta)$.

 Now  we apply the above mentioned  Gromov Lemma [2.4.1.A, Gr1986]
to get a  oriented curve $S^1 (t)$ in the linear space $T_x M$ such that
$$T({(\p /\p t) S^1(t)\over  | (\p /\p t) S^1 (t)|} ) = \rho + \eps\leqno (4.5.2) $$
for all $t$.
Next we observe that   for all $ \alpha > 0$ the curve $\alpha\cdot S^1(t)$ has
the same norm as $S^1 (t)$, i.e.
$$\Mm^1 _x ( T_{|  (\alpha \cdot S^1)}(t)) = \Mm^1 _x(T_{|  (S^1)}(t)) = \rho + \eps.$$
Thus we can assume  that our curve $S^1(t)$, which satisfies (4.5.2),
lies in the ball $B(0, \delta)$. By our choice of $\delta$ ( see (4.5.1)), we get
from (4.5.2)
$$\rho +{3\over 4} \eps \le \Mm^1 (Exp (S ^1(t)) )\le \rho + {5\over 4} \eps ,\leqno(4.5.3)$$
for all $t$.
This curve  $Exp (S^1(t))$ is an immersed curve. \QED


\medskip

Now let  us  to continue the proof of Proposition 4.4.a. We denote by $S^1(t)$
the embedded curve in Lemma 4.5. Next by choosing a tubular neighborhood
of $S^1(t)$ we can get a (small, thin)  oriented embedded solid torus  $T^3(t,s,r) =
S^1 (t) \times S^1 (s)\times [0,R]$
in $ M^m$ such that   our embedded curve is exactly the mean curve
$S^1(t)\times \{ 0\}\times \{ 0\}$
on the  solid torus. We can choose this torus $T^3$ so thin, such that for all $s,t,r$
we have
$$\Mm ^1 ( T^2_r(t,s)) \ge \rho + {\eps \over 4}.\leqno(4.5.4)$$
Using (4.5.4) we choose a smooth  unit vector field  $V(t,s)$ on the
torus $T^3(t,s, r)$  which is  tangential to each
torus $T^2 _r (t,s)$ such that $ T (V, V, V ) = \rho$. The integral curve
of this vector field is  either a circle or an curve of infinite length.
If there exists an integral curve of infinite length, then this curve is our
desired curve for the Proposition 4.5. Assume now that all the integral
curves are circles.  Then there exist an embedding $ S^1(t)\times [0,\mu]\times  [0,\mu] $
such that for all $(s, r)\in [0,\mu]\times [0,\mu]$ the circle  $ S^1 (t)\times \{ s\}
\times \{ r\}$ is an integral
curve of $V$. Now we perturb $V$ in  a neighborhood $[0, \alpha]\times
[0,\mu] \times [0,\alpha]$
with a very small $\alpha$  such that the perturbed  unit vector field $V'$  satisfies
$T(V', V', V') =\rho$ and the integral curve of vector field $V'$ is not
any more periodic. This completes the proof of the first part in Proposition 4.4.

\medskip

Using the same argument we  can prove the  second part b) of Proposition 4.4.
First we get the existence of an embedded curve  $S^1(t)$ of arbitrary length
on  $M$
such that $\Mm^1 ( T_{| S^1 } (t)) \ge \rho + (1/4)\eps$. Now
we consider a  torus tubular neighborhood of this curve in $M$ and apply the same
argument in the first part, namely we get on each torus $T^2 (t,s)$
an integral curve whose unit tangential vector $V = (\p / \p t) S^1 (t; s, r)$
satisfies the condition:
$$ T(V,V,V) =\rho.$$
If there exists an infinite integral curve, then we are done. If not, that
means all integral curve are circles, then we apply the perturbation
method in the proof of the first part and get our desired curve. \QED

\medskip

\section{ Existence of  isostatistical embeddings into ${\bm Cap ^ N}$.}

{\bf Main Theorem.} {\it  Any  compact smooth ($C^1$ resp.) statistical manifold $(M, g, T)$
can be immersed into  the statistical manifold $(Cap^N_+, g^F, T^{A-C})$  for some
finite number $N$. Hence any statistical manifold is a statistical model.}

We first deduce  our Main Theorem  from Theorem 5.1 and Theorem 5.5.

\medskip

{\bf 5.1. Theorem.} {\it Let $(M^m, g, T)$ be a compact smooth ($C^1$ resp.) statistical manifold. Then there exist
 numbers $N \in \N^+$ and  $A \ge 0$  as well as a smooth ($C^1$ resp.) embedding $f:
(M^m, g, T) \to (\R^{N}, g_0,A\cdot T_0)$ such that 
$f^* (g_0) = g$ and $f^* (A\cdot T_0) = T$. }

\medskip

Our proof of Theorem 5.1 uses the Nash embedding theorem, the Gromov embedding
theorem and an algebraic trick. The existence of monotone invariants prevents
us extend Theorem 5.1 for non-compact case (in contrast to the Riemannian  case.)

\medskip

{\bf 5.2. The Nash embedding theorem.} [Nash1954, Nash1956]  {\it Any smooth  ($C^1$ resp.) -Riemannian manifold
$(M^n, g)$ can be isometrically embedded into  $(\R ^{N}, g_0)$
for some $N$ depending on $M^n$.}

\medskip

We denote by $T_0$ the ``standard" 3-tensor on $\R^n$:
$$ T_0 = \sum _{i =1 } ^n dx_i ^3.$$

\medskip

{\bf 5.3. The Gromov immersion theorem.} [Gromov1986, 2.4.9.3' and 3.1.4] {\it Suppose that
$M^m$ is given with a smooth ($C^1$ resp.) symmetric 3-form $T$. Then there exists
an embedding $f : M^m \to \R^ {N_1(m)}$ with $N_1(m) = 3 (n + (^{n+1}_2) + ( ^{n+2}_3))$ such that $f^* (T_0) = T$.}

\medskip

{\it Proof of Theorem 5.1.} First we shall take an immersion
$f_1 : (M^m, g, T)
\to (\R^ {N_1(m)}, g_0, T_0)$ such that
$$ f_1 ^* (T_0) =  T.$$
The existence of $f_1$ follows from the Gromov immersion theorem.

Then we choose a positive number $A^{-1}$ such that
$$g - A^{-1} ( f _1 ^* ( g_0)) = g_1$$
is a Riemannian metric on $M$, i.e. $g_1$ is a positive  symmetric bi-linear form. Such a number $A$ exists, since  $M$
is compact.

Now we shall choose an  isometric immersion $f_2 : (M^m, g_1) \to (\R^{N} , g_0)$.
The existence of $f_2$ follows from the Nash isometric immersion theorem.

\medskip

{\bf 5.4. Lemma.} {\it  
There is a linear isometric  embedding $L_{m+1} : \R^{m+1} \to  \R^{2m +2}$
such that $L_{ n+1} ( T_0) = 0$.}

\medskip

{\it Proof.} We put
$$ L_{m+1}( x_1, \cdots, x_{m+1}) = ( f ^ 1( x_1), \cdots , f ^{m+1} ( x_{m+1}))$$
where $ f ^ i $ embeds the line $ (\R , ( dx ^ i) ^2, 0)$ into $( \R ^ 2 , ( dx  ^ { 2i -1})  ^ 2 + (dx ^{2i} ) ^ 2,
(dx ^ { 2i -1} ) ^ 3 + ( dx ^{ 2i} ) ^ 3)$:
$$ f^i (x_i) = {1\over \sqrt 2} (x_{2i-1} - x_{2i}).$$
 Clearly, $L_{m+1}$  is the required embedding.\QED

\medskip

{\it Completion of the proof of Theorem 5.1.}
Finally we take an embedding
$$ f_3 : M^m \to \R^{ (m+1)(m+2)+ m}  \oplus \R^{2m+2}$$
as follows.
$$ f_3(x) = A^{-1}\cdot f_1 ( x) \oplus (L_{n+1} \circ f_2).$$

Since $f_2$ is an embedding, $f_3$ is the required embedding map for Theorem 5.1.\QED

\medskip

{\bf 5.5. Theorem.} {\it  Suppose that $C$ is a compact
subset in $Cap ^{4n}_+$. Then any  bounded domain  $D$ in a linear statistical manifold
$ (\R^n, g_0,  A\cdot T_0)$
can be realized as an immersed statistical submanifold of $(Cap ^{4n}_+, g^F, T^{A-C}) $.}

\medskip
 Set
$$ T ^*: = \sum_{i =1} ^{n } {2 dx_i ^3\over  x_i}.$$

We denote by $S^n_{r, +}$ the positive  sector of the sphere  of radius $r$  centered at the origin in $\R^{n+1}$.

\medskip

{\it Proof of Theorem 5.5}.    We  choose  a very large positive number 
\begin{equation}
\bar A = \bar A(n, A)  \label{eq:bara}
\end{equation}
to be specified in Lemma \ref{lem:imm} later. 
 First, $\bar A$ in (\ref{eq:bara}) is required  to be so large  such that there exists a
number $1<\lambda< 2$  satisfying   the following equation
\begin{equation}
\lambda^2 +\frac{ 3 n}{(2 \bar A)^{2}} = 4. \label{eq:3.4}
\end{equation}

Equation (\ref{eq:3.4}) implies that $(\lambda, ( 2 \bar A )^{-1}, (2 \bar A)^{-1} , (2 \bar A)^{-1})\in \R^4$  is  a point in  $ S ^3_{2/\sqrt n, +}$. Hence there exists  a  positive number $r( \bar A)$ such that   for all $0<r \le r(\bar A)$ the ball $U(\bar A, r)$  of radius $r$  in  the sphere $S^3_{2/\sqrt n}$ that is centered  at the
point $(\lambda, (2 \bar A )^{-1}, (2 \bar A ) ^{-1},
(2 \bar A ) ^{-1})$  belongs  also to the  positive quadrant $ S^3_{2/\sqrt n, +}$. Hence $U(\bar A, r) \times _{ \text{ n  times } }  U(\bar A, r)$ is a subset in
$S_{2, +} ^{4n-1} \subset \R^{4n}$.

Next, we note that  Theorem 5.5 is a consequence of the following

\begin{lemma}\label{lem:imm}  For given positive numbers $R >0$ and $A\ge 0$ there exist  a   positive number $\bar A$, satisfying (\ref{eq:3.4}) and  depending only on $n$ and $A$,
 a positive number $r<r(\bar A)$   and  an isostatistical immersion $h$ from the bounded domain $[0,R]\times _{\text{ n times}} [0,R]\subset (\R^n, g_0, A\cdot T_0)$
into  $(Cap_+^{4N}, g^F, T^{A-C})$  such that  $h: [0,R]\times _{\text{ n times}} [0,R] \subset U(\bar A, r) \times _{ \text{ n  times } }  U(\bar A, r)$.
\end{lemma}

\begin{proof}  Set
$$ T ^*: = \sum_{i =1} ^{n } {2 dx_i ^3\over  x_i}.$$
Since $(U(\bar A, r), (g_0)|_{U(\bar A, r)}, T^*|_ {U(\bar A, r)})$ is a statistical submanifold of
$(\R^4_+, g_0,  T^*)$,  the direct product
$$(U(\bar A, r) \times _{\text{n  times}}U(\bar A, r), \oplus_{i =1} ^n (g_0)|_{U(\bar A, r)}, \oplus _{i =1} ^n T^*|_ {U(\bar A, r)})$$
 is a statistical submanifold of $(\R^{4n}_+, g_0,
T^*)$.  Since $(Cap_+^N, g^F, T^{A-C})$ is a statistical submanifold  of $(R^N_+, g_0, T^*)$ ,  we conclude that
$$(U(\bar A,r)\times _{\text{ n  times} } U(\bar A, r), \oplus_{i =1} ^n (g_0)|_{U(\bar A, r)}, \oplus _{i =1} ^n  T^*|_ {U(\bar A, r)})$$
 is a statistical submanifold of
$(Cap_+^{4N}, g^F, T^{A-C})$.  Hence,  to prove Lemma \ref{lem:imm}, it suffices to show that there are   positive numbers $\bar A = \bar A(n, A)$, $r< r(\bar A)$  and an isostatistical immersion
$f: ([0, R], dx ^2, A\cdot dx ^3)\to (U(\bar A, r), (g_0)|_{U(\bar A, r)}, T^*|_{U(\bar A, r)})$. On $U (\bar A, r)$
we consider the distribution $D(\rho)$  defined by
$$ D_x (\rho)  : = \{ v \in T_x U(\bar A, r): |v|_{g_0} = 1, T^* (v,v, v) = \rho\} $$
for  any given $\rho > 0$.
Clearly the existence of an isostatistical immersion $f: ([0, \R], dx ^2, A\cdot dx ^3)\to( U(\bar A, r), (g_0)|_{U(\bar A, r)}, T^*|_{U(A, r)})$ is equivalent to the existence of an integral curve with the
length $R$
of the distribution $D(A)$ on  $U (\bar A, r)$.  

Now  we  are going to  prove the following

\begin{lemma}\label{lem:torus} There exist  a positive number  $\bar A = \bar A(n, A)$  and an embedded torus $T^2$ in $U(\bar A, r)$ which
is provided with a unit vector field $V$ on $T^2$ such that $T ^*(V,V,V) = A$.
\end{lemma}

\begin{proof}[Proof of Lemma \ref{lem:torus}] Let us denote 
$$x_0 : = (\lambda, (2\bar A) ^{-1}, (2 \bar A) ^{-1},(2 \bar A ) ^{-1})\in S^3(2/\sqrt n)$$
 with $\lambda$ defined by (\ref{eq:3.4}). We shall need the following

\begin{lemma} \label{lem:sub2} There  exists  a  positive  number $\bar A = \bar A(n, A)$  such that the following    assertion holds. Let $H$ be any 2-dimensional subspace in $T_{x_0} U(\bar A, r) \subset \R^4$.
Then there exists a unit vector $w\in H$ such that $T^* (w,w,w) \ge \sqrt 2 A$.
\end{lemma}

\begin{proof}[Proof of Lemma \ref{lem:sub2}] Denote  by $\vec x_0$ the  vector in $\R^4$ with the same   coordinates as  those of  the point $x_0$.
For any given  $H$  as in Lemma \ref{lem:sub2}  there  exists   a unit vector $\vec h$  in $\R^4$,
which is not co-linear with $\vec x_0$ and which is orthogonal to $H$, such that
 a vector $w\in \R^4$  belongs to $H$ if and only if $w$  is a solution to following two linear equations:
\begin{equation}
 \la w, \vec x_0 \ra = 0 ,\label{eq:3.8}
 \end{equation}
 \begin{equation}
\la w, \vec h\ra  = 0. \label{eq:3.9}
\end{equation}
Adding a multiple of $\vec x_0$   to $\vec 
h$ if necessary, and  taking the  normalization, we  can  assume that
$$ \vec h = ( 0=h_1, h_2, h_3, h_4) \text { and } \sum_i h_i^2 = 1.$$
\medskip

\underline { Case 1}.  Suppose that not all the coordinates  $h _ i$  of $\vec h$ are of the same sign, so w.l.o.g.
we assume that $h_1 = 0, h_2  \le 0, h_3 > 0$. We put
\begin{eqnarray}
 k_2 : = {- h_2 \over \sqrt {(h_2) ^2+(h_3) ^2}}, \; k_3 : = {h_3 \over \sqrt{(h_2)^2 + 
(h_3) ^2}},\nonumber\\
 w : = ( w_1, w_2= (1 -\eps _2) k_3,w_3= (1-\eps _2)k_2,0=w_4)\in \R^4.\label{eq:3.10}
\end{eqnarray}
Obviously,  for any  choice of $w_1$ and $\eps_2$ the  equation (\ref{eq:3.9}) for $w$ is  satisfied. Now we choose  $w _1, \eps _2$ to be  solutions of the following  equations
\begin{equation}
 \lambda \cdot w _1 + ( 1-\eps _2)\cdot (2\bar A ) ^{-1}\cdot ( k_2 + k_3) = 0, \label{eq:3.11}
 \end{equation}
\begin{equation}
 w _1^ 2  = (2\eps _2 - \eps _ 2 ^ 2 ).\label{eq:3.12}
 \end{equation}
Note that (\ref{eq:3.11}) is equivalent to (\ref{eq:3.8}) and  (\ref{eq:3.12})  normalizes $w$.
From (\ref{eq:3.11})   we get
\begin{equation}
w _1 =-{(1 - \eps _2) ( k_2 + k_3) \over \lambda \cdot 2\bar A}.\label{eq:3.13}
\end{equation}
Substituting the value of $w_1$  into (\ref{eq:3.12}), we get
$$( { ( k_2 + k_3) ^2\over (\lambda \cdot 2\bar A )^2} + 1) \eps _2 ^ 2 -
( 2 + { 2(k_2 + k_3)^2\over (\lambda \cdot 2\bar A )^2 } ) \eps _2  + ({k_2 + k_3\over \lambda \cdot 2\bar A  }) ^2 = 0,$$
which we  simplify as follows:
\begin{equation}
\eps_2^2  -2 \eps_2 + \frac{(k_2 +k_3)^2 }{(k_2 + k_3) ^2 + 4 \lambda ^2 \bar A ^ 2}=0.\label{eq:3.14}
\end{equation}
Clearly, the following   choice of $\eps_2$  is a solution to (\ref{eq:3.14})
\begin{equation}
\eps _2 = 1-\frac{2\lambda \bar A}{\sqrt{(k_2 + k_3) ^2 + 4 \lambda ^2 \bar A ^2 }}.\label{eq:3.15}
\end{equation}
By our assumption on $h_2$ and $h_2$,  we have $0\le k_2, k_3 \le 1$. Since $1 < \lambda < 2$  by (\ref{eq:3.4}), we conclude that  when  $\bar  A$  goes  to infinity, the value $\eps_2$ goes to  zero.  Hence  there exists a number $N_1 >0$ such that  if $\bar A > N_1$ then
\begin{equation}
 \eps_2 > 0 \text{ and } (1-\eps_2) ^2 \ge {3\over 4}.\label{eq:eps2a1}
\end{equation}
  We shall show that  for $\eps_2 $ in (\ref{eq:3.15}) that  also  satisfies  (\ref{eq:eps2a1}) if $\bar A$ is  sufficiently large,  and for $w_1$ defined by (\ref{eq:3.13}),  the vector $w$ defined by (\ref{eq:3.10}) satisfies  the required condition   of Lemma \ref{lem:sub2}). Since $x_0 = ( \lambda, (2\bar A)^{-1},  (2\bar A)^{-1}, (2\bar A)^{-1})$ we have
\begin{equation}
 T^*_{x_0} (w, w, w) = {2 w_1 ^3\over \lambda} + (4\bar A) (w_2 ^3 + w_3 ^3).\label{eq:tw}
 \end{equation}
 Now assume that $\bar A >N_1$.   Noting that $\eps_2$ is positive and   close to zero,   and using $k_2 \ge 0$, $k _3 \ge 0$, we obtain from  (\ref{eq:3.10})
\begin{equation}
w_2 \ge 0,\,  w_3 \ge   0.\label{eq:est31}
\end{equation}
Since $0 < \eps < 1$,   $ 0 < k_2 + k_3 < 2$,  and $\lambda,  \bar A$  are positive,   we  obtain   from  (\ref{eq:3.13}) 
\begin{equation}
 w_1 < 0 \text{ and } |w_1| < \frac{1 }{\lambda \bar A}.\label{eq:est32}
\end{equation} 
Taking into account  (\ref{eq:3.10})  and (\ref{eq:eps2a1}), we  obtain
 \begin{equation}
   w_2 ^2 + w_3 ^2 = (1-\eps_2) ^2 \ge {3\over 4}.\label{eq:est33}
 \end{equation}
Using   (\ref{eq:est32}),  we obtain  from (\ref{eq:tw})
\begin{equation}
T^*_{x_0} (w, w, w)\ge {-2\over \lambda ^4 \bar A ^3}  +  (4\bar A)\cdot (w_2 ^3  + w_ 3^ 3). \label{eq:bara11}
\end{equation}
Observing  that the function $x^{3/2} +  (c- x) ^{3/2}$  is convex on interval $[0, c]$ for any $c>0$,    using    (\ref{eq:est31})  and  (\ref{eq:est33}), 
we obtain  from (\ref{eq:bara11})  
\begin{equation}
 T^*_{x_0} (w, w, w)\ge  {-2\over \lambda ^4 \bar A ^3} +   (4 \bar A) \cdot 2 (\frac{\sqrt 3}{ \sqrt 2}) ^3) =  {-2\over \lambda ^4 \bar A ^3} + 8(\sqrt{\frac{ 3}{2 }})^3   \bar A  . \label{eq:bara1}
\end{equation}
Increasing  $\bar  A$ if necessary,  noting that $1 < \lambda = \lambda(A)$,  equation  (\ref{eq:bara1}) implies that there exists a  large   positive number  $\bar  A(n, A)$   
depending only on $n$ and $A$ such that  any  subspace $H$ defined by the equations  (\ref{eq:3.8}) and (\ref{eq:3.9}),  where $h$  is in  Case 1,
contains  a unit vector $w$ that  satisfies the condition in Lemma \ref{lem:sub2}, i.e. the RHS of (\ref{eq:bara1})   is larger than $\sqrt 2 A$.
\medskip

\underline{ Case 2}.  W.l.o.g. we  assume that
$h_2 \ge h_3 \ge h_4 > 0$  and  therefore  we have
\begin{equation}
\alpha : = {h_2 + h_3\over h_4} \ge 2.\label{eq:3.17a}
\end{equation} 

 We  shall search  the required       vector $w$ for Lemma \ref{lem:sub2}  in the  following form

\begin{equation}
 w : = (w_1,  w_2=- ( 1-\eps _2),w_3= -( 1-\eps _2),w_4= \alpha  (1-\eps_2)).\label{eq:3.16}
 \end{equation}
The equations (\ref{eq:3.16}) and (\ref{eq:3.17a})    ensure   that  $\la w,\vec 
h \ra = 0$  for any choice of  parameters $( w_1, \eps_2)$  of $w$ in (\ref{eq:3.16}).
Next  we   require that the parameters $(w_1,\eps _2)$ of $w$  satisfy the following two equations
\begin{equation}
  \lambda \cdot w_1  + {( 1-\eps _2) (\alpha -2)\over 2\bar A}= 0,\label{eq:3.18}
\end{equation}
\begin{equation}
 w _1 ^ 2  + (1-\eps_2)^2  (2+\alpha ^2 ) = 1 .\label{eq:3.19a}
\end{equation}
Note  that (\ref{eq:3.18}) is equivalent to (\ref{eq:3.8})  and (\ref{eq:3.19a}) normalizes   $w$. 
From (\ref{eq:3.18}) we   express $w_1$ in terms of $\eps_2$ as follows
\begin{equation}
 w _1 = -\frac{ ( 1 - \eps_2)(\alpha -2)}{ \lambda 2\bar A}.\label{eq:3.20}
 \end{equation}
 Set 
\begin{equation}
B: =(2+ \alpha ^2) + \frac{(\alpha -2)^2 }{ 4\lambda^2 \bar  A^2}. \label{eq:B}
\end{equation}
Plugging  (\ref{eq:3.20}) into (\ref{eq:3.19a})  and  using (\ref{eq:B}),
we  obtain the following  equation for $\eps_2$
$$(1-\eps_2) ^2  B  -1 = 0, $$
which is equivalent  to the following equation
\begin{equation}
(1-\eps_2) ^2     = \frac{1}{ B}.\label{eq:eps2a}
\end{equation}
Since  $\alpha \ge 2$  by (\ref{eq:3.17a}),  from (\ref{eq:B})  we have $B > 0$. Clearly
\begin{equation}
\eps_2 : = 1 - \frac{1}{ \sqrt B}\label{eq:eps2}
\end{equation}
is a solution  to (\ref{eq:eps2a}).

Since  $ \alpha  \ge 2$  and $ \eps_2 \le 1$  by (\ref{eq:eps2}), we obtain  from (\ref{eq:3.20})  that $w_1 \le 0$.    Taking into account  $1<\lambda$, $\bar A >0$,  we derive  from (\ref{eq:3.20}) and  (\ref{eq:eps2}) the following  estimates
$$ T^*_{x_0} (w,w,w)= \frac{ 2w_1 ^3}{\lambda} + (4\bar A)(1-\eps_2)^3(\alpha ^3-2) $$
$$ > 2 w_1  ^3 +  ( 4 \bar A) (\alpha^3 -2) (1-\eps_2) ^3  $$
$$ =  \frac{(\alpha -2) ^3}{ 4 \bar A ^3(\sqrt B) ^3}+ 4\bar A \frac{(\alpha  ^3 -2)}{(\sqrt{B})^3} $$
$$\ge - \frac{\alpha ^3 -2} {4\bar A ^3(\sqrt B) ^3}+ 4\bar A \frac{(\alpha  ^3 -2)}{(\sqrt{B})^3}\: (\text{  since } \alpha \ge 2)$$
\begin{equation}
=  \frac{\alpha ^3-2}{(\sqrt B) ^3}( - \frac{ 1}{ 4 \bar A ^3} + 4\bar A ).\label{eq:est1}
\end{equation}

\begin{lemma}\label{lem:est4} There  exists  a large  number $\bar  A = \bar A(n, A)$ depending   only on  $n$ such that  for all choice  of $\alpha\ge 2$  we have 
$$\frac{ (\alpha  ^3 -2)}{(\sqrt{B})^3} \ge {1\over 10^2}.$$
\end{lemma}
\begin{proof}
To  prove  Lemma \ref{lem:est4} it suffices  to show that  for $\alpha \ge 2$ we have
\begin{equation}
10^4 (\alpha ^3 -2) ^2  \ge  B^3.\label{eq:100}
\end{equation}
Clearly there  exists  a positive number $N_2$ such that   if $\bar A > N_2$, then by (\ref{eq:B}),  we have   
\begin{equation}
 B < {3\over 2} ( 2 + \alpha ^2) \label{eq:Best}
\end{equation}
for any  $\alpha \ge 2$. Hence  (\ref{eq:100})  is  a consequence  of the following  relation
\begin{equation}
10^4 (\alpha  ^3 -2)^2 \ge [{3\over 2}(2 +\alpha ^2 )]^ 3, \label{eq:32}
\end{equation}
which  we   shall establish  now.
To  prove  (\ref{eq:32})  it suffices  to show that
\begin{equation}
10 ^3  (\alpha ^3 -2 )^2 \ge (2 + \alpha ^2 ) ^3. \label{eq:1000}
\end{equation}
 The   inequality (\ref{eq:1000})  is equivalent  to the following
 \begin{equation}
 999 \alpha  ^6  - 6  \alpha ^4  - 4000 \alpha ^3 -  12\alpha ^2  + 3992 \ge 0.\label{eq:992}
 \end{equation}
 Since  $\alpha  \ge 2$ it follows  that  $\alpha ^3 \ge 8 $ and  hence  
\begin{equation}
999 \alpha  ^ 6 - 4000 \alpha ^3 =  499 \alpha ^ 6 +  500 \alpha ^3 (\alpha ^ 3 - 8) \ge 499 \alpha ^6.\label{eq:estn1}
\end{equation}
  Using  $2 \alpha ^6 \ge  6 \alpha^4$, we  obtain
 \begin{equation}
 499 \alpha  ^6  - 6 \alpha ^4  \ge  497   \alpha^6 .\label{eq:estn2}
\end{equation}
Using  $a^4 \ge 16$, we obtain
\begin{equation}
497  \alpha^6 - 12 \alpha   ^ 2  =  496 \alpha ^6 + \alpha^2 (\alpha ^4 -12) > 496 \alpha ^6.\label{eq:estn3} 
\end{equation}
From  (\ref{eq:estn1}), (\ref{eq:estn2}), (\ref{eq:estn3})  we obtain
\begin{equation}
 999 \alpha  ^6  - 6  \alpha ^4  - 4000 \alpha ^3 -  12 \alpha ^2  + 3992 \ge   496 \alpha ^ 6 + 3992 >0.
 \end{equation}
 This  proves (\ref{eq:1000})  and   hence  completes  the   proof of Lemma \ref{lem:est4}.
\end{proof}

Lemma \ref{lem:est4} implies that   when $\bar A= \bar A(A, n)$  is sufficiently large,  the  RHS of (\ref{eq:est1}) is larger than $\sqrt 2A$. 
This  proves  the existence of $\bar  A$, which depends only on $n $ and $A$, for Case  2.

This  completes  the proof  of  Lemma \ref{lem:sub2}.

\end{proof}

From Lemma \ref{lem:sub2}  we obtain immediately the following.

\begin{corollary}\label{cor:torus}  The exists  a small neighborhood $U_1\ni  x_0$ in $\bar  U(\bar A, r)$ such that  the following   statement  holds.
For any $x \in  U_1$   and any  two-dimensional  subspace $H \subset T_{x} U_1$  we have
$$ \max \{ T^* (v,v,v)|\,  v\in H \text{ and } |v|_{g_0} = 1\}\ge {5\over 4} A. $$
\end{corollary}

{\it Completion of the proof of Lemma \ref{lem:torus}.} Let $\bar A = \bar A(n, A)$      satisfy the condition of Lemma \ref{lem:sub2}.
Now we choose a small embedded
torus $T^2$ in $U_1 \subset U(\bar A,r)$.   By   Corollary \ref{cor:torus}, for all $x\in T^2$ we have
\begin{equation}
 \max \{ T^* (v,v,v)|\,  v\in T_x T^2 \text{ and } |v|_{g_0} = 1\}\ge {5\over 4} A.\label{eq:3.21}
 \end{equation}

Denote by $T_1T^2$ the bundle of the unit  tangent   vectors  of $T^2$. Since  $T^2 = \R^2 / \Z^2$ is parallelizable,  we have  $T_1 T^2 = T^2 \times S^1$. Thus
the  existence  of a  vector field $V$  required  in Lemma \ref{lem:torus} is equivalent  to  the existence of a function $T^2 \to  S^1$   satisfying
the condition  of   Lemma \ref{lem:torus}.   Next  we claim that  there  exists    a unit  vector field  $W$  on  $T^2$ such that
$T^* (W, W, W) = 0$.  First  we choose  some orientation for $T^2$, that  induces  an orientation on  $T_1T^2$  and hence on the  circle $S^1$.
Take an arbitrary unit vector  field $W'$ on $T^2$, equivalently we pick a function $W':T^2 \to S^1$.  Now we consider  the   fiber bundle $F$ over $T^2$  whose  fiber over $x \in T^2$
consists  of the  interval  $[W', - W'] $ defined  by the  chosen  orientation on the  circle  of    unit  vectors  in $T_x S^2$. Since $T^* (W',W',W') = -T^* (W,W,W)$, for  each $x\in T^2$  there  exists  a value
$W$ on $F(x)$  such that $T^* (W, W, W) = 0$ and $ W$ is  closest to $W'$.   Using $W$ we identify  the circle $S^1$ with  the interval $[0,1)$.
The  existence  of  $W$  implies that   the existence  of a
function $V:T^2 \to  [0,1)$,  regarded  as  a unit vector field $V$ on $T^2$,   that satisfies
the condition  of   Lemma \ref{lem:torus}  is  equivalent  to  the existence of a   function $f: T^2 \to  [0, 1)$   satisfying  the same  condition.  Now let $V(x)$   be   the smallest value of     unit vector $V(x) \in [0, 1)\subset  S^1 (T_xT^2)$     such that
$$T^*(V(x), V(x),  V(x)) = A$$
  for   each $x\in T^2$.  The existence of $V(x)$ follows  from (\ref{eq:3.21}).
This  completes  the  proof  of Lemma \ref{lem:torus}.
\end{proof}

As we have noted,  Lemma \ref{lem:torus} implies  Lemma  \ref{lem:imm}.
\end{proof}
This finishes  the proof  of Theorem 5.5.
\QED

\begin{proof}[Proof of  Main Theorem] The existence of  an isostatistical
immersion   of a  compact statistical manifold $(M, g, T)$ into $(Cap_+^N, g^F, T^{A-C})$ for   some finite  $N$ follows from  Theorem
5.1 and Theorem 5.5.

\end{proof}

\medskip



\begin{theorem}\label{thm:embedding} Any smooth ($C^1$ resp.)   compact statistical manifold $(M^n, g, T)$
admits an isostatistical embedding into  the statistical manifold $(\Pp_+([N]), g^F, T^{A-C})$  for some
finite number $N$. 
\end{theorem}

\begin{proof} 
To   prove  Theorem \ref{thm:embedding}
we repeat   the  proof of Main Theorem, replacing the Nash immersion theorem by  the Nash embedding theorem. First  we observe  that our immersion $f_3$
constructed in   the proof of Lemma 5.4 is an embedding, if $f_2$ is an isometric embedding.
The  existence   of an isometric  embedding  $f_2$ is ensured   by the Nash theorem.
 Hence, if $M^n$ is compact,  to prove  the existence   of  an isostatistical embedding of $(M^n, g, T)$  into $(\Pp_+([N]), g^F, T^{A-C})$ it suffices   to prove
 the strengthened   version of Theorem 5.5, where the existence of    an isostatistical immersion  is replaced by
the  existence  of  an isostatistical
embedding. 

Recall that  the proof   of Theorem 5.5  is  reduced to the  proof   of the existence  of  an  isostatistical immersion of a  bounded    statistical  interval $([0, R],   dt ^2,  A\cdot dt ^3)$  into     a torus  $T^2$  of a small domain in $(S^3_{2/\sqrt n, +}, g_0, T^*)$,  see the proof of Lemma  \ref{lem:imm}.   Here      for  simplicity of notation,  we abbreviate the restriction
of $T^*$ to the  sphere  in consideration as  $T^*$.

The  statistical immersion      produced  with the help of Lemma  \ref{lem:torus}   will  be an embedding  if  not all the  integral curves of the distribution  $D(A)$  on  the torus $T^2$
are  closed curves.    Now  we shall   search  for  an   isostatistical embedding
of $([0, R],   dt ^2,  A\cdot dt ^3)$  into     a torus  $T^2 \times  T^2$  of a small domain in $(S^3_{1/\sqrt n, +}, g_0, T^*)  \times  (S^3 _{1/\sqrt n, +}, g_0, T^*) \subset (\R^8, g_0,    T^*)$.
Since $T^4$ is  parallelizable,  repeating the argument   at the end of the  proof of Lemma \ref{lem:imm}, we  choose a distribution  $D(A) \subset TT^4$   such that
$D(A) =  T^4 \times  S^2$    and 
$$D_x A =  \{  v \in T_x T^4 |\,  | v|_{g_0}  = 1, \text{  and } T^* ( v, v, v) = A \}. $$
Now assume that   the  integral  curves  of  $D(A)$     that lie on the  first   factor  $T^2 \times  y$ for all $y \in S^3 _{1/\sqrt n, +}$  are closed.
Since  $T^2$  is  compact, there is a   positive number $p_1$ such that  the  periods  of  these  integral curves  are at least $p_1$.

 Now    let  us consider  the following  integral  curve  $\gamma (t)$ of  $D(A)$ on $T^4$.    The curve $\gamma (t)$ begins  at  a point  $(0, 0, 0, 0) \in  T^4$. Here  we identify  $T^1$ with
$[0,1] / ( 0 = 1) $.   The   integral curve  lies  on $T^2 \times (0, 0)$ until   it approaches  $(0, 0, 0, 0)$ again.  Since $D_x (A) = S^2$, we can slightly modify the   direction  of $\gamma (t)$
and   let  it    leave the torus  $T^2 \times (0, 0)$   and after   a  very short  time   $\gamma (t)$ must stay  on  the torus  $T^2 \times (\eps, \eps)$ where  $\eps$ is sufficiently small.
W.l.o.g.  we assume  that   the period of any     closed  curve  of  the  distribution $D(A)  \cap  T(T^2\times (\eps, \eps))$ is at least  $p_1$.   Repeating  this   procedure,  since $R$ and $p_1$ are finite,
we      produce  an  embedding  of $([0, R],   dt ^2,  A\cdot dt ^3)$  into      $T^4 \subset (S^3_{1/\sqrt n, +}, g_0, T^*)  \times  (S^3 _{1/\sqrt n, +}, g_0, T^*) $.

This completes  the  proof   of Theorem \ref{thm:embedding}. 


\end{proof}

\end{document}